\documentclass[12 pt,twoside]{article}
\usepackage{graphicx}

\newtheorem{theorem}{\bf Theorem}[section]

\newtheorem{lemma}[theorem]{\bf Lemma}
\newtheorem{corollary}[theorem]{\bf Corollary}

{{\footnotesize\sc }}
\input{amssym}

\begin{document}
\title{{\large Around a conjecture of Erd\H{o}s on graph Ramsey numbers}}
\vspace{1cm}\author{\small  L. Maherani$^{\footnotesize\textrm{a}}$, G.R. Omidi$^{\footnotesize\textrm{a},\textrm{b},\textrm{1}}$\medskip \\
\\
\footnotesize  $^{\textrm{a}}$Department of Mathematical Sciences, Isfahan
University of Technology,\\ \footnotesize Isfahan, 84156-83111, Iran\\
\footnotesize  $^{\textrm{b}}$School of Mathematics, Institute for
Research
in Fundamental Sciences (IPM),\\
\footnotesize  P.O.Box: 19395-5746, Tehran, Iran\\
\footnotesize{l.maherani@math.iut.ac.ir, romidi@cc.iut.ac.ir}}
\date {}

\footnotesize\maketitle\footnotetext[1] {\tt This research was in part supported
by a grant from IPM (No.90050049)} \vspace*{-.9cm} \footnotesize
\begin{abstract}\rm{}
For given graphs $G_1$ and $G_2$ the Ramsey number
$R(G_1, G_2)$, is the smallest positive integer $n$
such that each blue-red edge coloring of the complete graph $K_n$ contains
a blue copy of $G_1$ or a red copy of $G_2$.
In 1983, Erd\H{o}s conjectured that there is an absolute constant $c$ such that $R(G)= R(G,G)\leq 2^{c\sqrt{m}}$ for any graph $G$ with $m$ edges and no isolated vertices. Recently this conjecture was proved by B. Sudakov. In this note, using the Sudakov's ideas we give an extension of his result and some interesting corollaries.\\
\noindent
\\{{Keywords}:{ \footnotesize Ramsey number, Erd\H{o}s conjecture, Complete graph.\medskip}}
\end{abstract}
\bigskip
\section{\normalsize{Introduction}}
For given graphs $G_1$ and $G_2$ the {\it Ramsey number}
$R(G_1, G_2)$, is the smallest positive integer $n$
such that  each blue-red edge coloring of  the complete graph $K_n$ contains
a blue subgraph isomorphic to $G_1$ or a red  subgraph isomorphic to $G_2$. We denote $R(G,G)$ by $R(G)$.
The existence of such a positive integer is guaranteed by Ramsey's classical result \cite{Ramsey}.
Since 1970's, Ramsey theory has grown into one of the most active
areas of research within combinatorics, overlapping variously with graph theory, number theory, geometry and logic.
Probably the most complicated question in this field is the estimating the Ramsey number of complete graphs. A basic result of Erd\H{o}s and Szekeres \cite{szekeres} implies the following theorem.

\begin{theorem}{\rm(\cite{szekeres})}\label{szekeres}
For every positive integer $n$, $R(K_{n}) \leq 2^{2n}$.
\end{theorem}

Using probabilistic methods, Erd\H{o}s \cite{erdos} obtained a lower bound for $R(K_n)$.

\begin{theorem}{\rm(\cite{erdos})}\label{erdos}
For every positive integer $n >2$, $R(K_{n}) \geq 2^{\frac{n}{2}}$.
\end{theorem}

Over the last sixty years, there have been several improvements on these bounds. The best upper bound is obtained by Conlon \cite{conlon}.
A problem on Ramsey numbers of general graphs was posed by Erd\H{o}s and Graham \cite{graham}, who conjectured that among all graphs with $m= {n \choose 2}$ edges and no isolated vertices, the complete graph on $n$ vertices has the largest Ramsey number. Since the number of vertices in a complete graph with $m$ edges is a constant multiple of $\sqrt{m}$ motivated by Theorem \ref{szekeres} and Theorem \ref{erdos}, Erd\H{o}s conjectured \cite{conj erdos} that there is a constant $c$ such that for all graphs $G$ with $m$ edges and with no isolated vertices, $R(G) \leq 2^{c \sqrt{m}}$. The first result in this direction, proved by Alon, Krivelevich and Sudakov \cite{alon}, is that $R(G) \leq 2^{c \sqrt{m} \log m}$. They also proved this conjecture for bipartite graphs. In 2011, Sudakov \cite{sudakov} gave a short and intelligent proof for this conjecture. In fact, he proved the following theorem.

\begin{theorem}{\rm(\cite{sudakov})}\label{sudakov}
If $G$ is a graph with $m$ edges and with no isolated vertices, $R(G) \leq 2^{250 \sqrt{m}}$.
\end{theorem}

Our main results in this note are the following two theorems.
The first theorem is an extension of Theorem \ref{sudakov} and will be proved by the same arguments in \cite{sudakov}
and the second one can be obtained by the first main theorem and a result in \cite{alon} on the Ramsey number of a bounded maximum degree graph and a complete graph. We give their proofs in the last section.
Through this note, the notations $\ln x$ and $\log x$ are the logarithms in the natural base e and 2, respectively.

\begin{theorem}\label{main}
Let $G_{i}$, $i=1,2$, be a graph of order $n_i$ with no
isolated vertices and $m$ be a positive integer with $2^{ \frac{106 \sqrt{m}}{\log m}}\geq n-27 \sqrt{m}$, where $n = \max \{ n_1,n_2 \}$. Also suppose that for each $27 \leq \alpha \leq \frac{\log^{3}m}{8}$ we
have the following properties:
\begin{itemize}
\item[I.]
There exists   $U_i \subseteq V(G_i)$, $i=1,2$, so that $\vert U_i \vert \leq \alpha \sqrt{m}$  and $\Delta (G_i - U_i )\leq \frac{2\sqrt{m}}{\alpha},$
\item[II.]
There exists   $V_i\subseteq V(G_i)$, $i=1,2$, so that $\vert V_i \vert \leq m^{\frac{3}{2}}$  and
 $\max \{ R(G_1-V_1,G_2 ),R(G_1 ,G_2-V_2 ) \} \leq 2^{36\sqrt{m}}$.
\end{itemize}
Then $R(G_1 ,G_2 )\leq 2^{250\sqrt{m}}$.
\end{theorem}

The join of two graphs $G$ and $H$, denoted by $G+H$, is a graph with vertex set $V(G\cup H)$ and edge set $E(G\cup H)\cup\{uv|u\in G, v\in H\}$.

\begin{theorem} \label{main2}
Let $G_1 =K_{p}$,  $m\geq 27$, $p \leq 27\sqrt{m}+\frac{16\sqrt{m}}{\log^{3}m}$ and $G_2 = K_{l}+H$ where $l \leq 27\sqrt{m}$, $\Delta (H)\leq \frac{16\sqrt{m}}{\log^{3}m}$, $n(H)\leq 2^{\frac{106\sqrt{m}}{\log m}}$ and there is $S\subseteq V(H)$ with $\vert S\vert \leq \sqrt{m^{3}}-27\sqrt{m}$ such that $\Delta (H-S)< \frac{\log m}{4}$. Then $$R(G_1,G_2)\leq 2^{250\sqrt{m}}.$$
\end{theorem}

Immediately we can obtained some corollaries using Theorems \ref{main} and \ref{main2}.

\begin{corollary}\label{edges m1 m2}
Suppose that $G_i, i=1,2$, is a graph with $m_i$  edges and with no isolated vertices. If $m= \max \{ m_1 , m_2 \}$, then $R(G_1 , G_2 )\leq 2^{250 \sqrt{m}}$.
\end{corollary}

 \medskip
\noindent\textbf{Proof. }Since $G_i$, $i=1,2$, has at most $2m$ vertices, then by the Theorem \ref{szekeres}, $R(G_1,G_2) \leq R(K_{2m},K_{2m}) \leq 2^{4m}$. If $m\leq 60^{2}$, then $4m\leq 250\sqrt{m}$ and so $R(G_1,G_2)\leq 2^{250\sqrt{m}}$. So we may assume that $m\geq 60^{2}$. Obviously $2^{ \frac{106 \sqrt{m}}{\log m}}\geq 2m-27 \sqrt{m} \geq n-27 \sqrt{m}$, where $n=\max \{n_1,n_2 \}$ and  $n_i$, $i=1,2$, is the number of vertices of $G_i$. Again clearly for $27 \leq \alpha \leq \frac{\log^{3}m}{8}$, the maximum degree
of a graph obtained by deleting $\alpha \sqrt{m}$ vertices of large degree of $G_i$ is at most $\frac{2m}{\alpha \sqrt{m}}=\frac{2\sqrt{m}}{\alpha}$. On the other hand, $m^{\frac{3}{2}}\geq n$. Hence the assertion holds by Theorem \ref{main}.
$\hfill \blacksquare$

\begin{corollary}\label{n vertices}
Suppose $G_i$, $i=1,2$, is a graph of order $n_i$ and $n = \max \{ n_1,n_2 \}$. If there exists  a subset $U_{i}\subseteq V(G_i )$, $i=1,2$, with $\vert U_i \vert \leq 27 \sqrt[3]{n}$ such that  $\Delta (G_i - U_i)\leq 54\frac{\sqrt[3]{n}}{\log^{3}n}$, Then $R(G_1 ,G_2 ) \leq 2^{250 \sqrt[3]{n}}$.
\end{corollary}
\medskip
\noindent\textbf{Proof. }The assertion holds by putting $m = n^{\frac{2}{3}}$ in Theorem \ref{main}.
$\hfill \blacksquare$

\bigskip

Using Theorem \ref{main2},  we get the following corollaries.

\begin{corollary}
Let $G_1 =K_{p}$, $m\geq 27$, $p \leq 27\sqrt{m} + \frac{16\sqrt{m}}{\log^{3}m}$ and $G_2 = K_{l}+q K_1$, where $l \leq 27\sqrt{m}$ and $q \leq 2^{106 \frac{\sqrt{m}}{\log m}}$.  Then $R(G_1 ,G_2)\leq 2^{250\sqrt{m}}$.
\end{corollary}

\begin{corollary}
For positive integers $p$ and $q$ with $q \leq 2^{\frac{53p}{27 \log \frac{p}{27}}}$, we have $R(K_{p},K_{p, q})\leq R(K_{p},K_{p}+q K_1)\leq 2^{\frac{250}{27}p}$.
\end{corollary}

\bigskip
\section{\normalsize Preliminaries}
In this section, we present some lemmas which will be used in the proof of Theorem \ref{main}.
Let $G=(V,E)$ be a graph and $U \subseteq V$. The induced subgraph of $G$ on $U$ is denoted by $G[U]$. The edge density of $G[U]$ is denoted by $d(U)$ and is defined by
$d(U)= \frac{e(U)}{{\vert U\vert \choose 2}},$
where $e(U)$ is the number of edges of $G[U]$. In an edge-coloring of $K_n$, we call an ordered pair $(X,Y)$ of disjoint subsets of vertices monochromatic if all edges in $X\cup Y$ incident to each vertex in $X$ have the same color.

\begin{lemma}{\rm(\cite{sudakov})}\label{base mono pair}
For all $k$ and $l$, every blue-red edge coloring of $K_{N}$, contains a monochromatic pair $(X,Y)$ with
$$\vert Y\vert \geq {k+l \choose k}^{-1}N-k-l,$$
which is red and $\vert X\vert =k$ or is blue and  $\vert X\vert =l$.
\end{lemma}

\begin{lemma}{\rm(\cite{sudakov})}\label{mono pair}
Let $0<\epsilon \leq \frac{1}{7}$ and let $t$ and $N$ be positive integers satisfying $t\geq \epsilon^{-1}$ and $N\geq t \epsilon^{-14\epsilon t}$. Then every blue-red edge coloring of $G=K_{N}$ in which  edge density of  the induced subgraph on the blue edges is at most $\epsilon$ contains a monochromatic pair $(X,Y)$ with $\vert X\vert \geq t$ and $\vert Y\vert \geq \epsilon^{14\epsilon t}N$.
\end{lemma}

\begin{lemma}{\rm(\cite{sudakov})}\label{sparse S}
Let $G$ be a graph with $n$ vertices, maximum degree $\Delta$ and let $\epsilon \leq \frac{1}{8}$. If $H$ has $N\geq  \epsilon^{4\Delta \log \epsilon }n$ vertices and does not contain a copy of $G$, then it has a subset $S$ with $\vert S\vert \geq \epsilon^{-4\Delta \log \epsilon }N$ and edge density $d(S)\leq \epsilon$.
\end{lemma}

\begin{lemma}\label{mono pair'}
Let $G_{1}$ and $G_2$ be graphs with no isolated vertices of orders $n_1$ and $n_2$, respectively, and $m$ be a positive integer such that $2^{ \frac{106 \sqrt{m}}{\log m}}\geq n-27 \sqrt{m}$, where $n = \max \{ n_1,n_2 \}$. Also suppose that for $27 \leq \alpha \leq \frac{\log^{3}m}{8}$
there exists  $U_i \subseteq V(G_i)$ with $\vert U_i \vert \leq \alpha \sqrt{m}$  such that   $\Delta (G'_i )\leq \frac{2\sqrt{m}}{\alpha}$, where $G'_i =G_i - U_i$.
If a blue-red edge coloring of $K_{N}$ has no blue copy of $G_1$ and no red copy of $G_2$ and it contains a monochromatic pair $(X,Y)$ with $\vert X\vert \geq \alpha \sqrt{m}$ and $\vert Y\vert \geq 2^{125\alpha ^{\frac{-1}{3}}\sqrt{m}}$, then it also contains a monochromatic pair $(X',Y')$ with $\vert X'\vert \geq 2^{2\alpha^{ \frac{1}{3}}}\sqrt{m}$ and $\vert Y'\vert \geq 2^{-120\alpha ^{\frac{-1}{3}}\sqrt{m}} \vert Y\vert$.
\end{lemma}
\medskip
\noindent\textbf{Proof. }Consider a blue-red edge coloring of $K_N$ with the properties in Lemma \ref{mono pair'}. Assume that the color of the  monochromatic pair $(X,Y)$ is blue. By the hypothesis, the induced subgraph  on $Y$ contains no blue $G'_1$. Otherwise, together with $X$, we get a blue $G_1$. Let $\epsilon = 2^{-3 \alpha^{\frac{1}{3}}}$  and $t=2^{2\alpha^{\frac{1}{3}}}\sqrt{m}$.  Since $\alpha \leq \frac{\log^{3}m}{8}$ we have $2^{\alpha^{\frac{1}{3}}} \leq \sqrt{m}$ and so $t\geq \epsilon^{-1}$. Also note that, since $27 \leq \alpha \leq \frac{\log^{3}m}{8}$, we have $42\alpha^{\frac{1}{3}} 2^{-\alpha^{\frac{1}{3}}} \leq 48\alpha^{\frac{-1}{3}}$ and $2^{5\alpha^{\frac{-1}{3}}\sqrt{m}} \geq 2^{10\frac{\sqrt{m}}{\log m}} \geq m^{\frac{3}{2}} \geq 2^{2\alpha^{\frac{1}{3}}}\sqrt{m} = t$.
Applying Lemma \ref{sparse S} to the blue graph restricted to $Y$, we find a subset $S\subset Y$ with
\begin{eqnarray*}
\vert S\vert &\geq & \epsilon^{-4\Delta (G'_1) \log \epsilon } \vert Y\vert \geq (2^{-3\alpha^{\frac{1}{3}}})^{-4(2\alpha^{-1}\sqrt{m})(-3\alpha^{\frac{1}{3}})}\vert Y\vert \\
 &= & 2^{-72\alpha^{\frac{-1}{3}}\sqrt{m}} \vert Y\vert \geq 2^{53\alpha^{\frac{-1}{3}}\sqrt{m}} \geq 2^{ \frac{106 \sqrt{m}}{\log m}}\geq n-27 \sqrt{m},
\end{eqnarray*}
such that the density of the induced subgraph on the blue edges in $S$ is at most $\epsilon$. Then the size of $S$ satisfies
$$\vert S\vert \geq 2^{53\alpha^{\frac{-1}{3}}\sqrt{m}} \geq 2^{5\alpha^{\frac{-1}{3}}\sqrt{m}} 2^{48\alpha^{\frac{-1}{3}}\sqrt{m}} \geq t2^{42\alpha^{\frac{1}{3}}2^{-\alpha^{\frac{1}{3}}}\sqrt{m}} = t\epsilon^{-14\epsilon t},$$
and we can apply Lemma \ref{mono pair} to $S$. So $S$ contains a monochromatic pair $(X',Y')$ with $\vert X'\vert \geq t$ and $\vert Y'\vert \geq \epsilon^{14\epsilon t} \vert S\vert$. To complete the proof, recall that $\vert S\vert \geq 2^{-72\alpha^{\frac{-1}{3}}\sqrt{m}} \vert Y\vert$ and therefore
$$\vert Y'\vert \geq \epsilon^{14\epsilon t} \vert S\vert \geq 2^{-48\alpha^{\frac{-1}{3}\sqrt{m}}} \vert S\vert \geq 2^{-120\alpha ^{\frac{-1}{3}}\sqrt{m}} \vert Y\vert.$$
A similar argument, which we omit, can be used to finish the proof in the case when the color of  the monochromatic pair $(X,Y)$ is red.
$\hfill \blacksquare$
\bigskip

\section{\normalsize The Proofs}
In this section, we give the proofs for the main theorems.\\

\noindent\textbf{ The proof of Theorem \ref{main}. }Assume that $N=250\sqrt{m}$ and suppose for contradiction that there is a blue-red edge coloring of $K_{N}$ with no blue copy of $G_1$ and no red copy of $G_2$. Since $G_1$ and $G_2$ have at most $n$ vertices by Theorem \ref{szekeres}, we have $R(G_1,G_2)\leq R(K_{n})\leq 2^{2n}$. So we may assume that $n\geq 125 \sqrt{m}$. Applying Lemma \ref{base mono pair} with $k=l= 27\sqrt{m}$, we have a monochromatic pair $(X_1,Y_1)$ with $\vert X_1\vert \geq 27\sqrt{m}$ and
$$\vert Y_1\vert \geq {k+l \choose k}^{-1}N-k-l \geq 4^{-27\sqrt{m}}N =2^{196\sqrt{m}}.$$
Define $\alpha_1=27$ and $\alpha_{i+1} = 2^{2\alpha_{i}^{\frac{1}{3}}}$. It is easy to see $\alpha_{i+1} \geq (\frac{4}{3})^{3}\alpha_i$ for all $i$ and therefore $\alpha_{i}^{-\frac{1}{3}} \leq \frac{1}{3}(\frac{3}{4})^{i-1}$. This implies that $$\Sigma_{j=1}^{i} \alpha_{j}^{-\frac{1}{3}} \leq \frac{1}{3}\Sigma_{j=0}^{i-1}(\frac{3}{4})^{j}\leq \frac{1}{3} \Sigma_{j\geq 0}(\frac{3}{4})^{j} =\frac{4}{3}.$$
Since the blue-red edge coloring of $K_{N}$ has no blue copy of $G_1$ and no red copy of $G_2$, we can repeatedly apply Lemma \ref{mono pair'}. After $i$ iterations, we have a monochromatic pair $(X_{i+1},Y_{i+1})$ with $\vert X_{i+1}\vert \geq \alpha_{i+1}\sqrt{m}$ and $$ \vert Y_{i+1}\vert \geq 2^{-120\alpha_{i} ^{\frac{-1}{3}}\sqrt{m}} \vert Y_{i}\vert \geq 2^{-120 \sqrt{m} \Sigma_{j=1}^{i} \alpha_{j}^{-\frac{1}{3}}} \vert Y_1 \vert \geq 2^{-120\sqrt{m}(\frac{4}{3})} 2^{196\sqrt{m}} = 2^{36\sqrt{m}}.$$
We continue iterations until the first index $i$ such that $\alpha_{i} \geq \frac{\log^{3}m}{8}$. Then for $\alpha = \frac{\log^{3}m}{8}$ we have a monochromatic pair $(X,Y)$, $\vert X\vert \geq \alpha \sqrt{m}$ and $\vert Y\vert \geq 2^{36\sqrt{m}} \geq 2^{125\alpha ^{\frac{-1}{3}}\sqrt{m}}$. Then applying Lemma \ref{mono pair'} one more time we obtain a monochromatic pair $(X',Y')$  with $\vert X'\vert \geq 2^{2\alpha^{\frac{1}{3}}}\sqrt{m} = m^{\frac{3}{2}}$ and $\vert Y'\vert \geq 2^{36\sqrt{m}}$. First let $(X',Y')$ be blue. By the second property in the hypothesis, since the induced subgraph on $Y'$ contains no a red $G_2$, it contains a blue copy of $G_1-V_1$ and clearly $X' \cup (G_1-V_1)$ contains a blue copy of $G_1$, a contradiction. So we may assume that the monochromatic pair $(X',Y')$ is red. Clearly the subgraph on $Y'$ does not contain a blue $G_1$ and so  it contains a red $G_2-V_2$ and so $X' \cup (G_2-V_2)$ contains a red copy of $G_2$, a contradiction.
$\hfill \blacksquare$\\

To prove Theorem \ref{main2} we need the following result in \cite{alon}.

\begin{theorem} (\cite{alon}) \label{alon}
Let $H$ be a graph with $h$ vertices and chromatic number $\chi(H)=k\geq 2$. Suppose that
there is a proper $k$-coloring of $H$ in which the degrees of all vertices, besides possibly those
in the first color class, are at most $r$, where $1\leq r (<h)$. Define $\alpha(k, r)$ to be 1 if $k>r$, and
0 otherwise. Then, for every integer $m>1$,
$$R(H,K_m)\leq(\frac{100m}{\ln m})^{\frac{(2r-k+2)(k-1)}{2}} (\ln m)^{\alpha(k, r)}h^{r}.$$
\end{theorem}

\bigskip

\noindent\textbf{The proof of Theorem \ref{main2}. }Let $q=n(H-S)$, $k=\chi(H-S)$ and $r=\Delta (H-S )$.
Since $m\geq 27$, then $p\leq 40\sqrt{m}$. If $r=0$, we add an edge, say e, to $H-S$.  By Theorem \ref{alon} we have
\begin{eqnarray*}
R(H-S,G_1) &\leq & R(H-S +e ,G_1) \\
&<&100 pq \leq 100 (27\sqrt{m}+\frac{16\sqrt{m}}{\log^{3}m})2^{106\frac{\sqrt{m}}{\log m}}\\
&<& 4000\sqrt{m} 2^{\frac{106\sqrt{m}}{\log m}}=2^{\log 4000\sqrt{m}+106\frac{\sqrt{m}}{\log m}} \\
&\leq & 2^{13+\frac{1}{2}\log{m} +106\frac{\sqrt{m}}{\log m}}< 2^{36\sqrt{m}}.
\end{eqnarray*}
Then we can use Theorem \ref{main}. Now suppose that $r\geq 1$. Using Theorem \ref{alon}, we get

\begin{eqnarray*}
R(H-S,G_1) = R(H-S,K_{p})&\leq &
(\frac{100 p}{\ln p})^{\frac{(2r-k+2)(k-1)}{2}} (\ln p) (q)^{r}\\
&\leq &(\frac{100 p}{\ln p})^{\frac{r^2 +r}{2}} (\ln p) (27\sqrt{m}+2^{106\frac{\sqrt{m}}{\log m}}-\sqrt{m^{3}})^{r} \\
&\leq &\frac{(100p)^{{\frac{r^2 +r}{2}}}2^{106\frac{\sqrt{m}}{\log m}r}}{(\ln
p)^{\frac{r^2 +r-2}{2}}}\\
&\leq & (4000\sqrt{m})^{\frac{r^2 +r}{2}}2^{106\frac{\sqrt{m}}{\log m}r}y\\
&=& 2^{(\frac{r^{2}+r}{2}\log 4000\sqrt{m} + 106\frac{\sqrt{m}}{\log m}r)}y \\
 &<&2^{36^{\sqrt {m}}}.
\end{eqnarray*}
The last inequality holds, since
 $$y=\frac{1}{(\ln p)^{\frac{r^2 +r-2}{2}}} \leq 1,$$
and
\begin{eqnarray*}
\frac{r^{2}+r}{2}\log 4000\sqrt{m} + 106\frac{\sqrt{m}}{\log m}r &=& \frac{r^2 +r}{2}(\log 4000 + \log \sqrt{m})+106\frac{\sqrt{m}}{\log m}r\\
& <& 6(r^2 +r) + \frac{r^{2}+r}{2} \log \sqrt{m} + 106\frac{\sqrt{m}}{\log m}r \\
&\leq &12(\log ^2 \sqrt[4]{m}) + \log ^2 \sqrt[4]{m} \log \sqrt {m} + 26.5 \sqrt{m}\\
&= &\frac{12}{16} \log ^2 m +\frac{1}{32}\log ^3 m+ 26.5 \sqrt{m} \\
&< & 9.5 \sqrt{m} +26.5 \sqrt{m} = 36\sqrt{m}.
\end{eqnarray*}

Now we use Theorem \ref{main} to finish the proof.
$\hfill \blacksquare$


\end{document}